\documentclass[a4paper]{article}
\usepackage{amsmath,amsthm,amssymb}
\usepackage{graphicx}
\newtheorem{Def}{Definition}
\newtheorem{theorem}[Def]{Theorem}
\newtheorem{lemma}[Def]{Lemma}
\newtheorem{corollary}[Def]{Corollary}
\newtheorem{remark}[Def]{Remark}
\newtheorem*{maintheorem}{Main Theorem}
\newtheorem*{problem}{Problem}
\makeatletter
 
\renewcommand{\@cite}[2]{{#1\if@tempswa , #2\fi}}
\makeatother

\title{Generating the mapping class group of a punctured surface by involutions}
\author{Naoyuki Monden}
\date{July 5, 2008}
\begin{document}
\maketitle
\begin{abstract}
Let $\Sigma_{g,b}$ denote a closed orientable surface of genus $g$ with $b$ punctures and let $\rm Mod(\Sigma_{\textit{g,b}})$ denote its mapping class group.
In [Luo] Luo proved that if the genus is at least 3, $\rm Mod(\Sigma_{\textit{g,b}})$ is generated by involutions. He also asked if there exists a universal
upper bound, independent of genus and the number of punctures, for the number of torsion elements/involutions needed to generate $\rm Mod(\Sigma_{\textit{g,b}})$.
Brendle and Farb [BF] gave an answer in the case of $g\geq 3, b=0$ and $g\geq 4, b=1$, by describing a generating set consisting of 6 involutions. Kassabov
showed that for every $b$ $\rm Mod(\Sigma_{\textit{g,b}})$ can be generated by 4 involutions if $g\geq 8$, 5 involutions if $g\geq 6$ and 6 involutions
if $g\geq 4$. We proved that for every $b$ $\rm Mod(\Sigma_{\textit{g,b}})$ can be generated by 4 involutions if $g\geq 7$ and 5 involutions if $g\geq 5$.
\end{abstract}

\section{Introduction}
\ \ \ Let $\Sigma_{g,b}$ be an closed orientable surface of genus $g\geq 1$ with arbitrarily chosen $b$ points (which we call punctures). Let $\rm Mod(\Sigma_{\textit{g,b}})$ be the mapping class group of $\Sigma_{g,b}$, which is the group of homotopy classes of
orientation-preserving homeomorphisms preserving the set of punctures. Let $\rm Mod^{\pm}(\Sigma_{\textit{g,b}})$ be the extended mapping class group of $\Sigma_{g,b}$, which is the group of homotopy class of all
(including orientation-reversing) homeomorphisms preserving the set of punctures. By $\rm Mod_{\textit{g,b}}^{0}$ we will denote  the subgroup of $\rm Mod_{\textit{g,b}}$ which fixes the punctures pointwise. It is clear that we have
the exact sequence:
\begin{center}
$1 \rightarrow  \rm Mod_{\textit{g,b}}^{0} \rightarrow \rm Mod_{\textit{g,b}} \rightarrow Sym_{\textit{b}} \rightarrow 1,$
\end{center}
where the last projection is given by the restriction of a homeomorphism to its action on the puncture points.\\
\ \ \ The study of the generators for the mapping class group of a closed surface was first considered by Dehn. He proved in [De] that $\rm Mod(\Sigma_{\textit{g},0})$ is generated by a finite set of Dehn twists.
Thirty years later, Lickorish [Li] showed that $3\textit{g}-1$ Dehn twists generate $\rm Mod_{\textit{g},0}$. This number was improved to $2\textit{g}+1$ by Humphries [Hu]. Humphries proved, moreover, that in fact the number $2\textit{g}+1$ is minimal;
i.e. $\rm Mod(\Sigma_{\textit{g},0})$ cannot be generated by $2g$ (or less) Dehn twists. Johnson [Jo] proved that the $2g+1$ Dehn twists also generate $\rm Mod(\Sigma_{\textit{g},1})$.
In the case of multiple punctures the mapping class group can be generated by $2g+b$ Dehn twists for $b\geq 1$ (see [Ge]).\\
\ \ \ It is possible to obtain smaller generating sets of $\rm Mod(\Sigma_{\textit{g,b}})$ by using elements other than twists. N.Lu (see [Lu]) constructed a generated  set of
$\rm Mod(\Sigma_{\textit{g},0})$ consisting of 3 elements. This result was improved by Wajnryb who found the smallest possible generating set of $\rm Mod(\Sigma_{\textit{g},0})$
consisting of 2 elements (see [Wa]). Korkmaz proved in [Ko] that one of these generators can be taken as a Dehn twist. It is also known that in the case of $b=0$ the mapping class group
can be generated by 3 torsion elements (see [BF]). More, Korkmaz showed in [Ko] that the mapping class group can be generated by 2 tosion elements (also in the case of $b=0,1$).
In [Ma], Maclachlan proved that the moduli space is simply connected as a topological space by showing that $\rm Mod(\Sigma_{\textit{g},0})$ is generated by torsion elements.
Several years later Patterson generalized these results to $\rm Mod(\Sigma_{\textit{g,b}})$ for $g\geq 3, b\geq 1$ (see [Pa]).\\
\ \ \ In [MP], McCarthy and Papadopoulos proved that $\rm Mod(\Sigma_{\textit{g},0})$ is generated by infinitely
many conjugetes of a single involution for $g\geq 3$. Luo, see [Luo], described the finite set of involutions which generate $\rm Mod(\Sigma_{\textit{g,b}})$ for $g\geq 3$. 
He also proved that $\rm Mod(\Sigma_{\textit{g,b}})$ is generated by torsion elements in all cases except $g=2$ and $b=5k+4$, but this group is not generated by involutions if $g\leq 2$.
Brendle and Farb proved that $\rm Mod(\Sigma_{\textit{g,b}})$ can be generated by 6involutions for $g\geq 3,b=0$ and $g\geq 4,b\leq 1$ (see [BF]). In [Ka], Kassabov proved that for every $b$
$\rm Mod(\Sigma_{\textit{g,b}})$ can be generated by 4 involutions if $g\geq 8$, 5 involutions if $g\geq 6$ and 6 involutions if $g\geq 4$. He also proved in the case of $\rm Mod^{\pm}(\Sigma_{\textit{g,b}})$.\\
Our main result is stronger than [Ka].
\begin{maintheorem}
For all $g\geq 3$ and $b\geq 0$, the mapping class group $\rm Mod({\Sigma_{\textit{g,b}}})$ can be generated by:\\
\ \ \ $(a)$ $4$ involutions if $g\geq 7$;\\
\ \ \ $(b)$ $5$ involutions if $g\geq 5$.
\end{maintheorem}

\section{Preliminaries}
\ \ \ Let $c$ be a simple closed curve on $\Sigma_{g,b}$. 
Then the (right hand) Dehn twist $T_{c}$ about $c$ is the homotopy class of the homeomorphism obtained
by cutting $\Sigma_{\textit{g,p}}$ along $c$, twisting one of the side by $360^{\circ}$ to the right and gluing two sides of a back to each ohter. Figure 1 shows the Dehn twist about the curve $c$. 
\begin{figure}[htbp]
 \begin{center}
  \includegraphics*[width=10cm]{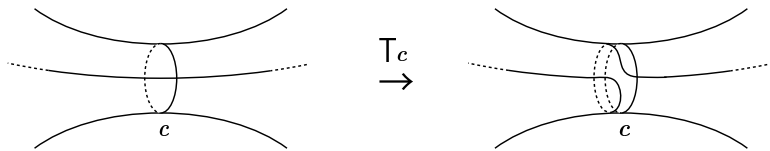}
 \end{center}
 \caption{The Dehn twist}
 \label{fig:one}
\end{figure}
We will denote by $T_{c}$ the Dehn twist around the curve $c$. \\
\ \ \ We record the following lemmas.
\begin{lemma}
For any homeomorphism $h$ of the surface
$\Sigma_{g,b}$ the twists around the curves $c$ and $h(c)$ are conjugate in the mapping class group $\rm Mod(\Sigma_{\textit{g,b}})$,
\begin{eqnarray*}
T_{h(\textit{c})}&=&hT_{\textit{c}}h^{-1}.
\end{eqnarray*}
\end{lemma}
\begin{lemma}
Let $c$ and $d$ be two simple closed curves on $\Sigma_{g,b}$. If $c$ is disjoint from $d$, then 
\begin{eqnarray*}
T_{c}T_{d}=T_{d}T_{c}
\end{eqnarray*}
\end{lemma}

\section{Proof of main theorem}
In this section we proof maintheorem. The keypoints of proof are to generate $T_{\alpha}$ in 4 involutions by using lantern relation.

\subsection{The policy of proof}
We give the policy of proof of maintheorem.
\begin{lemma}
Let $G$, $Q$ denote the groups and let $N, H$ denote the subgroups of $G$. Suppose that the group $G$ has the following exact sequence;
\begin{eqnarray*}
1 \rightarrow N \xrightarrow[]{i} G \xrightarrow[]{\pi} Q \rightarrow 1.
\end{eqnarray*}
If $H$ contains $i(N)$ and has a surjection to $Q$ then we have that $H=G$.
\end{lemma}
\begin{proof}
We suppose that there exists some $g\in G-H$. By the existence of surjection from $H$ to $Q$, we can see that there exists some $h\in H$ such that $\pi(h)=\pi(g)$.
Therefore, since $\pi(g^{-1}h)=\pi(g)^{-1}\pi(h)=1$, we can see that $g^{-1}h\in \rm Ker \ \pi= Im \ \textit{i}$. Then there exists some $n\in N$ such that $i(n)=g^{-1}h$.
By $i(N)\subset H$, since $i(n)\in H$ and $h\in H$, we have
\begin{eqnarray*}
g=h\cdot i(n)^{-1}\in H.
\end{eqnarray*}
This is contradiction in $g\notin H$. Therefore, we can prove that $H=G$.
\end{proof}
It is clear that we have the exact sequence:
\begin{eqnarray*}
1 \rightarrow  \rm Mod_{\textit{g,b}}^{0} \rightarrow \rm Mod_{\textit{g,b}} \rightarrow Sym_{\textit{b}} \rightarrow 1.
\end{eqnarray*}
Therefore, we can see the following corollary;
\begin{corollary}
Let $H$ denote the subgroup of $\rm Mod(\Sigma_{\textit{g,b}})$, which contains $\rm Mod^{0}(\Sigma_{\textit{g,b}})$ and has a surjection to $\rm Sym_{\textit{b}}$. Then H is equal to $\rm Mod(\Sigma_{\textit{g,b}})$.
\end{corollary}
We generate the subgroup $H$ which has the condition of corollary 4 by involutions.\\
\ \ \ \ Let us embed our surface $\Sigma_{g,b}$ in the Euclidian space in two different ways as shown on Figure 2. (In these pictures we will assume that
genus $g=2k+1$ is odd and the number of punctures $b=2l+1$ is odd. In the case of even genus we only have to swap the top parts of the pictures, and in the case of even
number of punctures we have to remove the last point.)\\
\ \ \ In Figure 2 we have also marked the puncture points as $x_{1},\ldots,x_{b}$ and we have the curves $\alpha_{i}$, $\beta_{i}$, $\gamma_{i}$ and $\delta$.
The curve $\alpha_{i}$, $\beta_{i}$, $\gamma_{i}$ are non separating curve and $\delta$ is separating curve.\\
\ \ \ Each embedding gives a natural involution of the surface---the half turn rotation around its axis of symmetry. Let us call these involutions $\rho_{1}$
and $\rho_{2}$.\\
\begin{figure}[htbp]
 \begin{center}
  \includegraphics*[width=13cm]{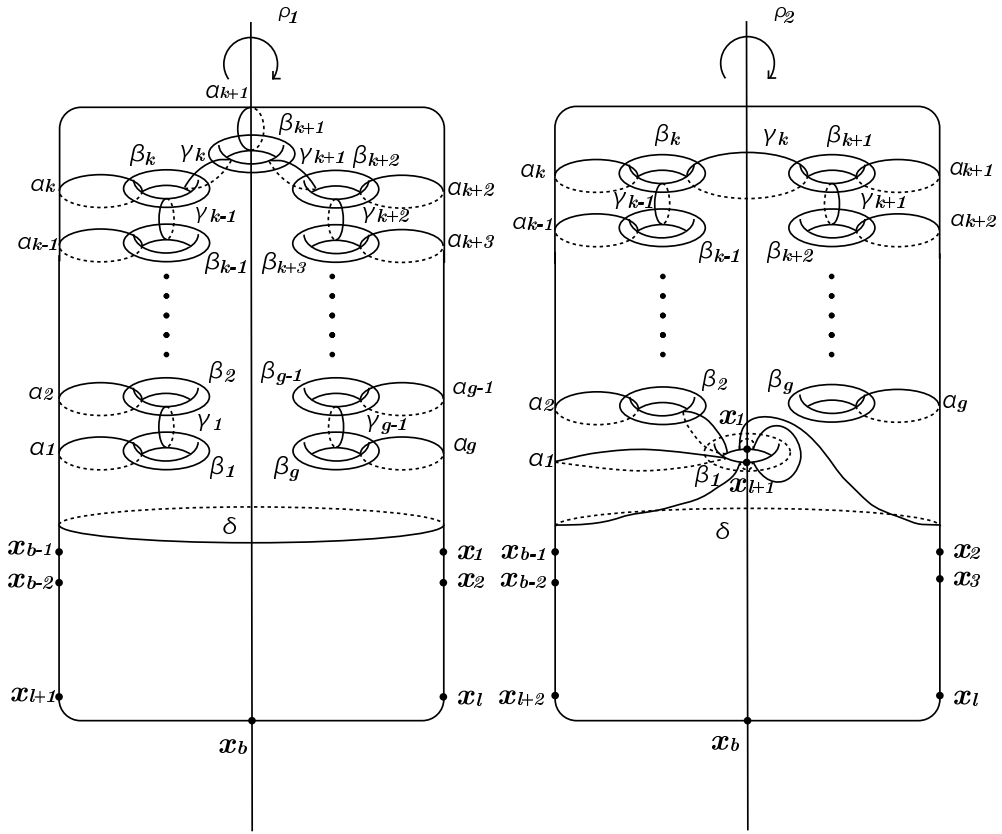}
 \end{center}
 \caption{The embeddings of the surface $\Sigma_{g,b}$ in Euclidian space used to define the involutions $\rho_{1}$ and $\rho_{2}$.}
 \label{fig:one}
\end{figure}
\ \ \ Then we can get following lemma;
\begin{lemma}
The subgroup of the mapping class group be generated by $\rho_{1}$, $\rho_{2}$ and 3 Dehn twists $T_{\alpha}$, $T_{\beta}$ and $T_{\gamma}$ around one
of the curve in each family contains the subgroup $\rm Mod^{0}(\Sigma_{\textit{g,b}})$.
\end{lemma}
We postopne the proof of lemma 5 until Section 5.\\
\ \ \ The existence a surjection from the subgroup $H$ of $\rm Mod(\Sigma_{\textit{g,b}})$ to $\rm Sym_{\textit{b}}$ is equivalent to showing taht the $\rm Sym_{\textit{b}}$
can be generated by involutions;
\begin{eqnarray*}
r_{1}&=&(1,\textit{b}-1)(2,\textit{b}-2)\cdots (\textit{l},\textit{l}+1)(\textit{b}) \\
r_{2}&=&(2,\textit{b}-1)(3,\textit{b}-2)\cdots (\textit{l},\textit{l}+2)(1)(\textit{l}+1)(\textit{b}) \\
r_{3}&=&(1,\textit{b})(2,\textit{b}-1)(3,\textit{b}-2)\cdots (\textit{l},\textit{l}+2)(\textit{l}+1)
\end{eqnarray*}
corresponding to 3 involutions in $H$.
\begin{lemma}
The symmetric group $\rm Sym_{\textit{b}}$ is generated by $r_{1}, r_{2}$ and $r_{3}$.
\end{lemma}
\begin{proof}
The group generated by $r_{i}$ contains the long cycle $r_{3}r_{1}=(1,2,\ldots,\textit{b})$ and transposition $r_{3}r_{2}=(1,\textit{b})$. These two elements generate the whole symmetric group, therefore
the involutions $r_{i}$ generate $\rm Sym_{\textit{b}}$.
\end{proof}
We note that the images of $\rho_{1}$ and $\rho_{2}$ to $\rm Sym_{\textit{b}}$ are $r_{1}$ and $r_{2}$.\\
\ \ \ Therefore, by Lemma 1, Corollary 4, Lemma 5 and Lemma 6 we sufficient to generate $H$ by $\rho_{1}$, $\rho_{2}$ and involutions which have the following conditions;
\begin{description}
\item[$\langle 1 \rangle$] involutions which genarate the Dehn twist around $\gamma$,
\item[$\langle 2 \rangle$] two of each involutions which exchange $\alpha$ and $\beta$, $\beta$ and $\gamma$, $\gamma$ and $\alpha$,
\item[$\langle 3 \rangle$] involution whose image is $r_{3}$.
\end{description}

\subsection{Generating Dehn twists by 4 involutions}
\ \ \ In this subsection, we argue about $\langle 1 \rangle$. Moreover, we generate Dehn twists by 4 involutions. The basic idea is to use the lantern relation.\\
\ \ \ We begin by recalling the lantern relation in the mapping class group. This relation was first discovered by Dehn and later 
rediscovered by Johnson. \\

\begin{figure}[htbp]
 \begin{center}
  \includegraphics*[width=4.5cm]{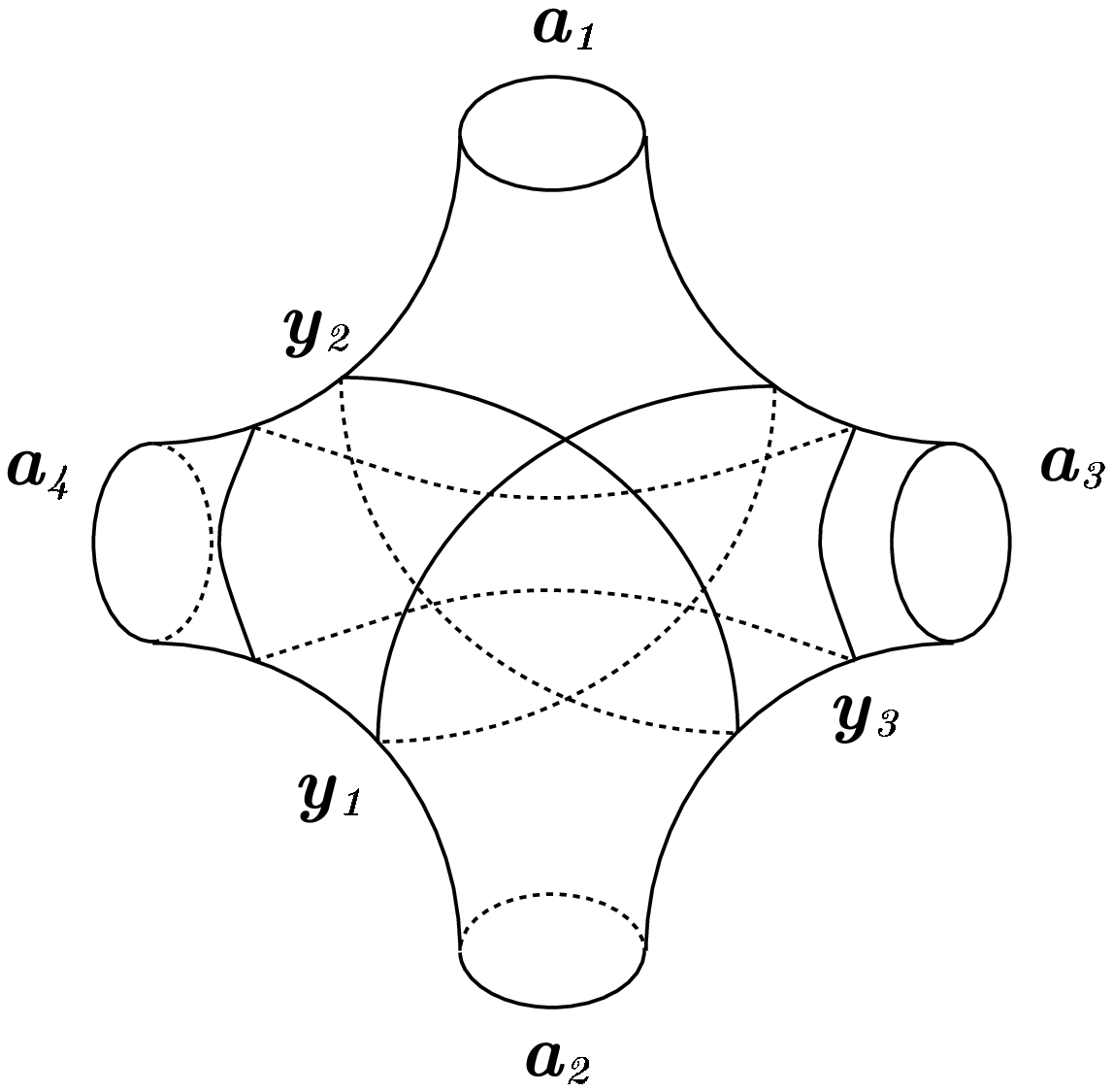}
 \end{center}
 \caption{Lantern}
 \label{fig:one}
\end{figure}
From now on we will assume that the genus $g$ of the surface is at least 5.\\
\ \ \ Let the $S_{0,4}$ be a surface of genus 0 with 4 boundary components. Denote by $a_{1}, a_{2}, a_{3}$ and $a_{4}$ the four boundary
curves of the surface $S_{0,4}$ and let the interior curves $y_{1}, y_{2}$ and $y_{3}$ be as shown in Figure 3.\\
\ \ \ The following relation:
\begin{eqnarray}
T_{y_{1}}T_{y_{2}}T_{y_{3}}=T_{a_{1}}T_{a_{2}}T_{a_{3}}T_{a_{4}}.
\end{eqnarray}
among the Dehn twists around the curves $a_{i}$ and $y_{i}$ is known as the lantern relation. Notice that the curves $a_{i}$ do not
intersect any other curve and that the Dehn twists $T_{a_{i}}$ commute with every twists in this relation. This allows us to
rewrite the lantern relation as follows
\begin{eqnarray}
T_{a_{4}}=(T_{y_{1}}T_{a_{1}}^{-1})(T_{y_{2}}T_{a_{2}}^{-1})(T_{y_{3}}T_{a_{3}}^{-1}).
\end{eqnarray}
\ \ \ Let $R$ denote the product $\rho_{2}\rho_{1}$. By Figure 2 we can see that $R=\rho_{2}\rho_{1}$ acts as follows:
\begin{eqnarray}
 R\alpha_{\textit{i}}&=&\alpha_{\textit{i}+1},  \ (1\leq \textit{i}<\textit{g})\nonumber \\
 R\beta_{\textit{i}}&=&\beta_{\textit{i}+1},   \ (1\leq \textit{i}<\textit{g})\\
 R\gamma_{\textit{i}}&=&\gamma_{\textit{i}+1},  \ (1\leq \textit{i}<\textit{g}-1).\nonumber
\end{eqnarray}
\ \ \ The lanterns $S$ and $R^{-2}S$ have a common boundary component $a_{1}=R^{-2}a_{2}$ and their union is a surface $S_{2}$ homeomorphic to a sphere with 6 boundary components. By Figure 4 we can see that there exists an involution $\bar{J}$ of $S_{2}$ which takes $S$ to $R^{-2}S$.\\
\ \ \ Let us embed the surface $S_{2}$ in $\Sigma_{g,b}$ as shown on Figure 5. The boundary components of $S_{2}$ are $a_{1}=\alpha_{k}$, $a_{2}=\alpha_{k+2}$, $a_{3}=\gamma_{k+1}$, $a_{4}=\gamma_{k}$, $R^{-2}a_{1}=\alpha_{k-2}$, $R^{-2}a_{2}=\alpha_{k}$, $R^{-2}a_{3}=\gamma_{k-1}$ and $R^{-2}a_{4}=\gamma_{k-2}$; and the middle curve $y_{1}=\alpha_{k+1}$.
The Figure 5 shows the existence of the involution $\tilde{J}$ on the complement of $S_{2}$ which is a surface of genus $g-5$ with 6 boundary components. Gluing together $\bar{J}$ and $\tilde{J}$ gives us the involution $J$ of the surface $\Sigma_{g,b}$. By Figure 4 $J$ acts as follows
\begin{eqnarray*}
J(a_{1})=R^{-2}a_{2}, \ J(a_{3})=R^{-2}a_{1}, \ J(y_{1})=R^{-2}y_{2}, \ J(y_{3})=R^{-2}y_{1}.
\end{eqnarray*}
Therefore, we have 
\begin{eqnarray}
R^{2}J(a_{1})=a_{2}, \ R^{2}J(y_{1})=y_{2} \nonumber \\
JR^{-2}(a_{1})=a_{3}, \ JR^{-2}(y_{1})=y_{3}.
\end{eqnarray}
Let $\rho_{3}$ denote $T_{a_{1}}\rho_{2}T^{-1}_{a_{1}}$.\\
By Lemma 1, (4) and that $\rho_{2}$ sends $a_{1}=\alpha_{k}$ to $y_{1}=\alpha_{k+1}$, we have
\begin{eqnarray}
T_{y_{1}}T_{a_{1}}^{-1}=\rho_{2}T_{a_{1}}\rho_{2}T^{-1}_{a_{1}}=\rho_{2}\rho_{3},\nonumber \\
T_{y_{2}}T^{-1}_{a_{2}}=R^{2}J\rho_{2}\rho_{3}JR^{-2}, \ T_{y_{3}}T^{-1}_{a_{3}}=JR^{-2}\rho_{2}\rho_{3}R^{2}J.
\end{eqnarray}
By (2) and (5) we have
\begin{eqnarray}
T_{\gamma_{\textit{k}}}=(\rho_{2}\rho_{3})(R^{2}J\rho_{2}\rho_{3}JR^{-2})(JR^{-2}\rho_{2}\rho_{3}R^{2}J).
\end{eqnarray}
\begin{figure}[htbp]
 \begin{center}
  \includegraphics*[width=7cm]{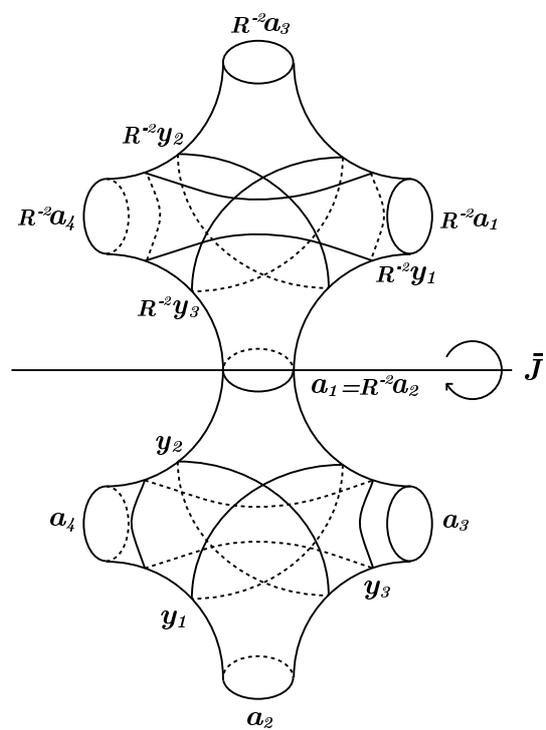}
 \end{center}
 \caption{$S_{2}$ and the involution $\bar{J}$}
 \label{fig:one}
\end{figure}
\begin{figure}[htbp]
 \begin{center}
  \includegraphics*[width=12cm]{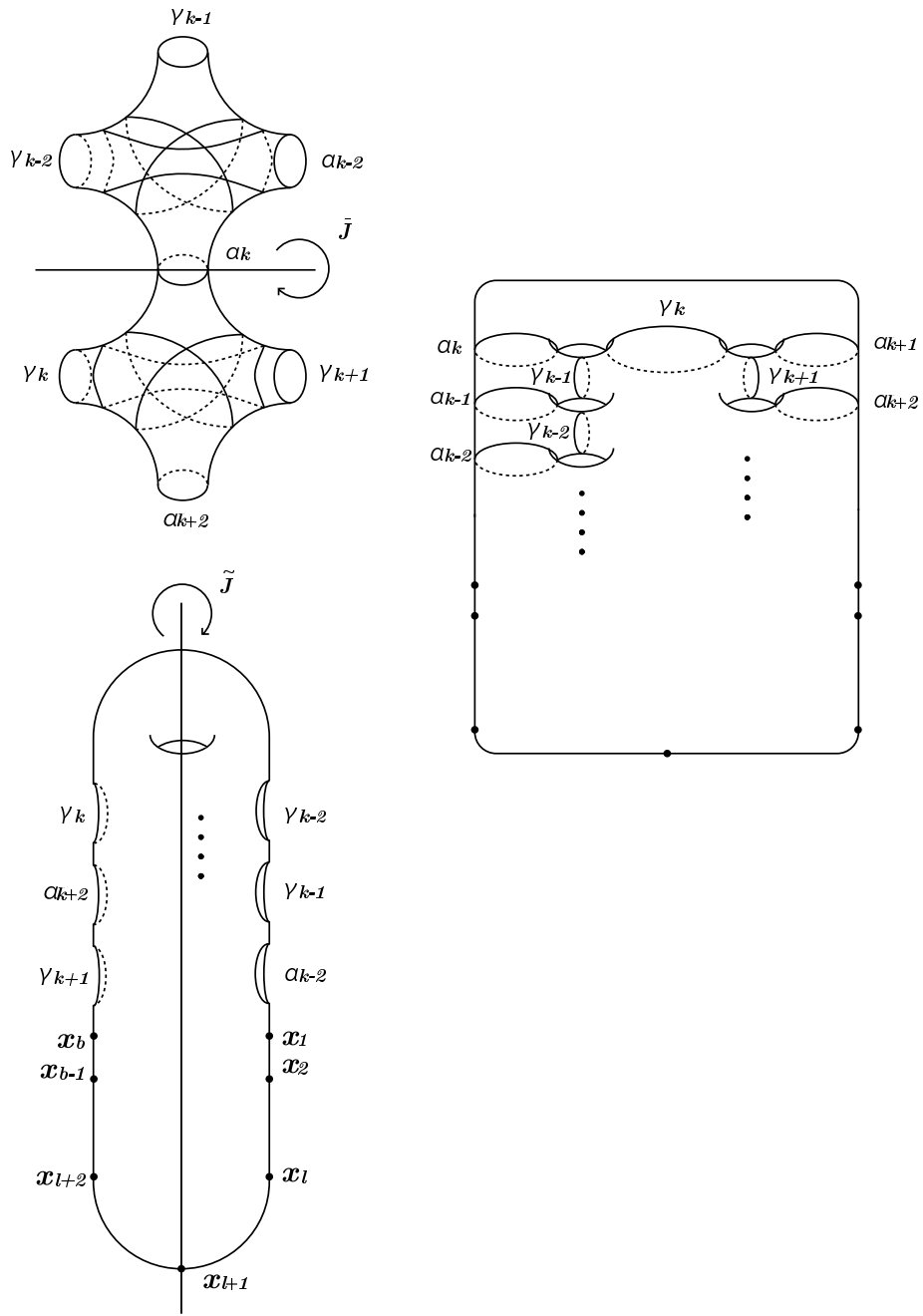}
 \end{center}
 \caption{The involution $J$ on $\Sigma_{g,b}$}
 \label{fig:one}
\end{figure}

\subsection{Genus at least 5}
\ \ \ We proof that the mapping class group is generated by 5 involutions.\\
\ \ \ The five involutions are $\rho_{1}, \rho_{2}, \rho_{3}, J$ and another involution $I$. We construct involution $I$ in the same way as involution $J$ like Figure 6.
\begin{figure}[h]
 \begin{center}
  \includegraphics*[width=8cm]{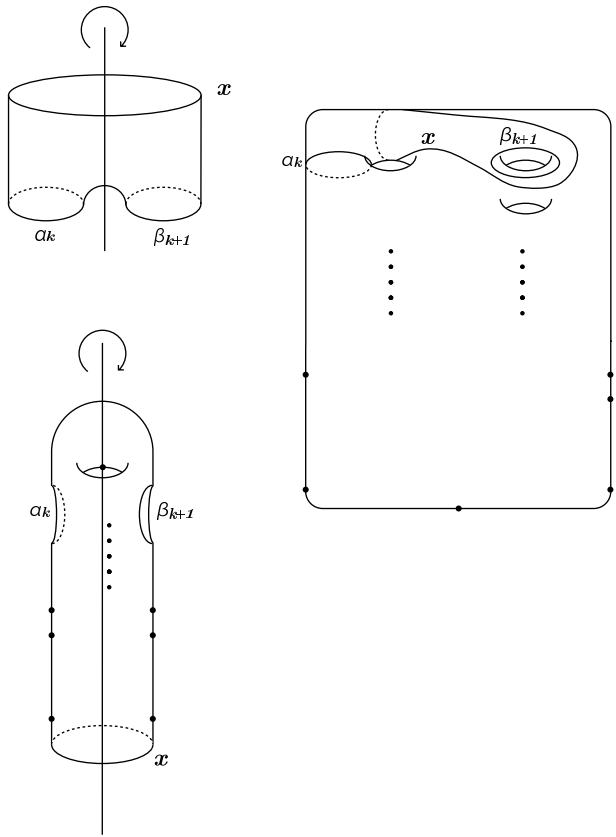}
 \end{center}
 \caption{The involution $I$ on $\Sigma_{g,b}$}
 \label{fig:one}
\end{figure}
\begin{theorem}
If $g\geq 5$, the group $G_{3}$ generated by $\rho_{1}$, $\rho_{2}$, $\rho_{3}$, $I$ and $J$ is the whole mapping class group $\rm Mod(\Sigma_{\textit{g,b}})$.
\end{theorem}
\begin{proof}
By the relation (6) we satisfy the condition $\langle 1 \rangle$. Since $J$ sends $\alpha_{k-2}$ to $\gamma_{k+1}$ and $I$ sends $\alpha_{k}$ to $\beta_{k+1}$, we consist the condition $\langle 2 \rangle$.
We can also see that we satisfy the condition $\langle 3 \rangle$ from a way to the construction of the involution $J$.\\
\ \ \ Therefore, we can finish the proof of the theorem because we can satisfy the conditions in 3.1.
\end{proof}

\subsection{Genus at least 7}

\begin{figure}[htbp]
 \begin{center}
  \includegraphics*[width=12cm]{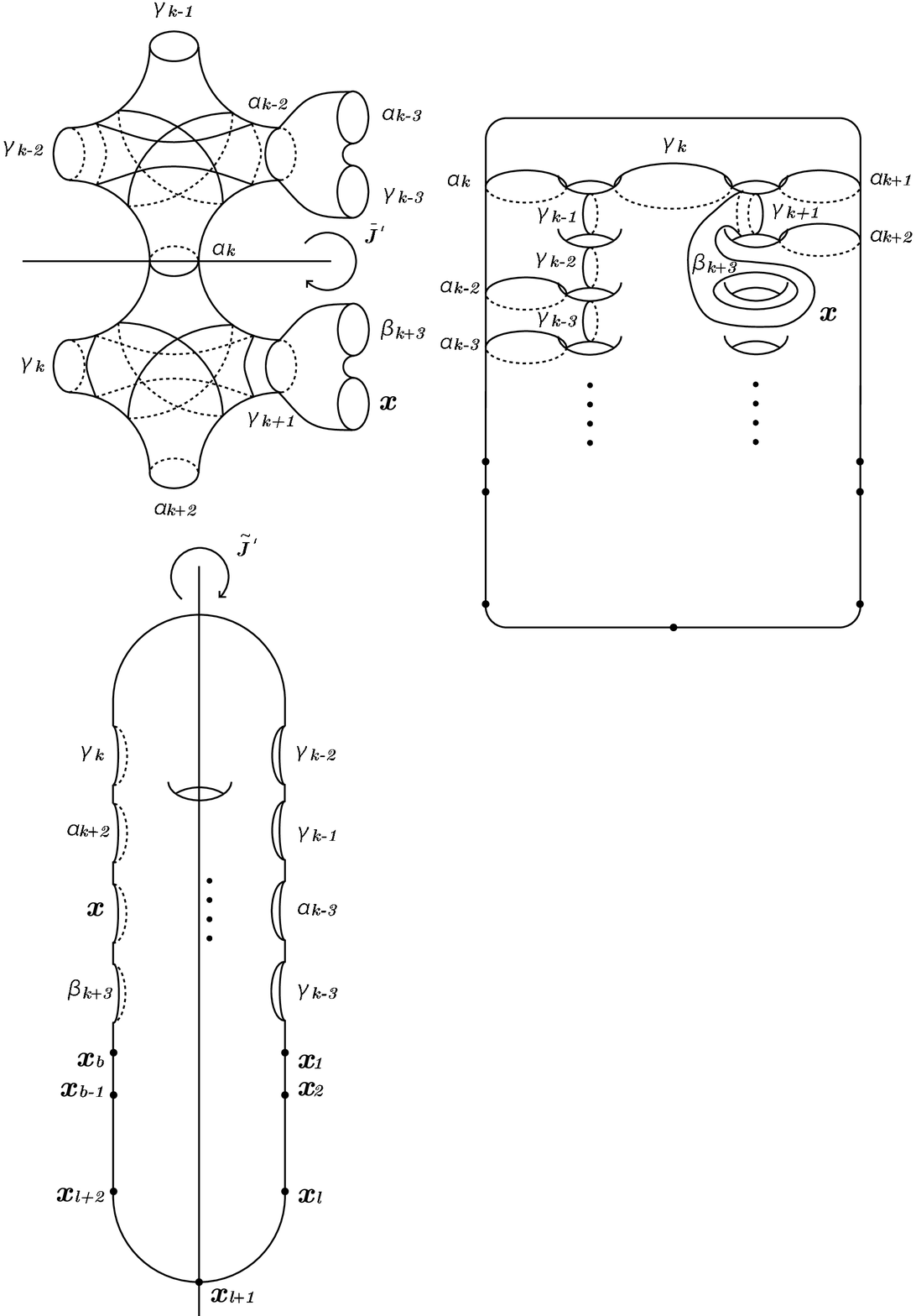}
 \end{center}
 \caption{The involution $J^{\prime}$ on $\Sigma_{g,b}$}
 \label{fig:one}
\end{figure}

\ \ \ We want to improve the above argument and show that for the genus $g\geq 7$ we do not need the involution $I$ in order to generate the mapping class group.
Assume that the genus of the surface is at least 7.\\
\ \ \ The $S_{2}$ and two pairs of pants have common boundary components $R^{-2}a_{1}$ and $a_{3}$ and their union is a surface $S_{3}$ homeomorphic to a sphere with 8 boundary components.
Figure 7 shows the existence of the involution $\bar{J^{\prime}}$ on $S_{3}$ which extends the involution $\bar{J}$ on $S_{2}$. \\
\ \ \ Let us embed $S_{3}$ in the $\Sigma_{g,b}$ as shown on Figure 7. From Figure 7 we can find the involution $\tilde{J^{\prime}}$ of the complement of $S_{3}$. Let $J^{\prime}$ be
the involution obtained by gluing together $\bar{J^{\prime}}$ and $\tilde{J^{\prime}}$. Moreover, from Figure 7 we can construct $J^{\prime}$ which acts on the punctures as the involution $r_{3}$. 

\begin{theorem}
If $g\geq 7$, the group $G_{4}$ generated by $\rho_{1}$, $\rho_{2}$, $\rho_{3}$ and $J^{\prime}$ is the whole mapping class group $\rm Mod(\Sigma_{\textit{g,b}})$.
\end{theorem}

\begin{proof}
From the construction of $J^{\prime}$ we have
\begin{eqnarray*}
T_{\gamma_{\textit{k}}}=(\rho_{2}\rho_{3})(R^{2}J^{\prime}\rho_{2}\rho_{3}J^{\prime}R^{-2})(J^{\prime}R^{-2}\rho_{2}\rho_{3}R^{2}J^{\prime}) \in G_{4}.
\end{eqnarray*}
Therefore, we can see that we satisfy the condition $\langle 1 \rangle$. Since $J^{\prime}$ can send $\alpha_{k-2}$ to $\gamma_{k+1}$ and $\beta_{k+3}$ to $\gamma_{k-3}$, we can satisfy the condition $\langle 2 \rangle$ only in $J^{\prime}$.
Moreover, By that $J^{\prime}$ acts as $r_{3}$, we consist the condition $\langle 3 \rangle$. Therefore, the group $G_{4}$ is the whole mapping class group.
\end{proof}

\section{The subgroup generated by 2 involutions and 3 Dehn twists, which contains $\rm Mod^{0}(\Sigma_{\textit{g,b}})$}
\ \ \ In this section we prove Lemma 5 (i.e. we construct the subgroup of the mapping class group $\rm Mod(\Sigma_{\textit{g,b}})$ by 2 involutions and 3 Dehn twists, which contains the pure mapping class group $\rm Mod^{0}(\Sigma_{\textit{g,b}})$).\\
\ \ \ We recall that $R=\rho_{2}\rho_{1}$. By Lemma 1 and (3), we get following relation;
\begin{eqnarray}
T_{\alpha_{\textit{i}+1}}&=&RT_{\alpha_{\textit{i}}}R^{-1} \nonumber \\
T_{\beta_{\textit{i}+1}}&=&RT_{\beta_{\textit{i}}}R^{-1}\\
T_{\gamma_{\textit{i}+1}}&=&RT_{\gamma_{\textit{i}}}R^{-1}. \nonumber
\end{eqnarray}
\ \ \ Let the subgroup $G$ of the mapping class group be generated by $\rho_{1}$, $\rho_{2}$ and 3 Dehn twists $T_{\alpha}$, $T_{\beta}$ and $T_{\gamma}$ around one
of the curve in each family. By relation (7), $T_{\alpha_{i}}, T_{\beta_{i}}, T_{\gamma_{i}} \in G$ for all $i$.\\
\ \ \ Our next step is to show that $G$ contains $\rm Mod^{0}(\Sigma_{\textit{g,b}})$.
Let denote the curves $\delta^{\prime},\eta^{\prime},\delta^{\prime\prime},\eta^{\prime\prime},\delta_{\textit{j}},\eta_{\textit{j}} (\textit{j}=1,\ldots,\textit{l}-1,\textit{l}+1,\ldots,\textit{b}-2)$ in Figure 8. 
In [Ge] it is shown that $\rm Mod^{0}(\Sigma_{\textit{g,b}})$ is generated by Dehn twists around the curves $\alpha_{i}$-es, $\beta_{i}$-es, $\gamma_{i}$-es, $\delta^{\prime}$, $\delta^{\prime\prime}$ and $\delta_{j}$-es, for $j=1,\ldots, l-1, l+1,\ldots ,b-2$.
\begin{figure}[htbp]
 \begin{center}
  \includegraphics*[width=13cm]{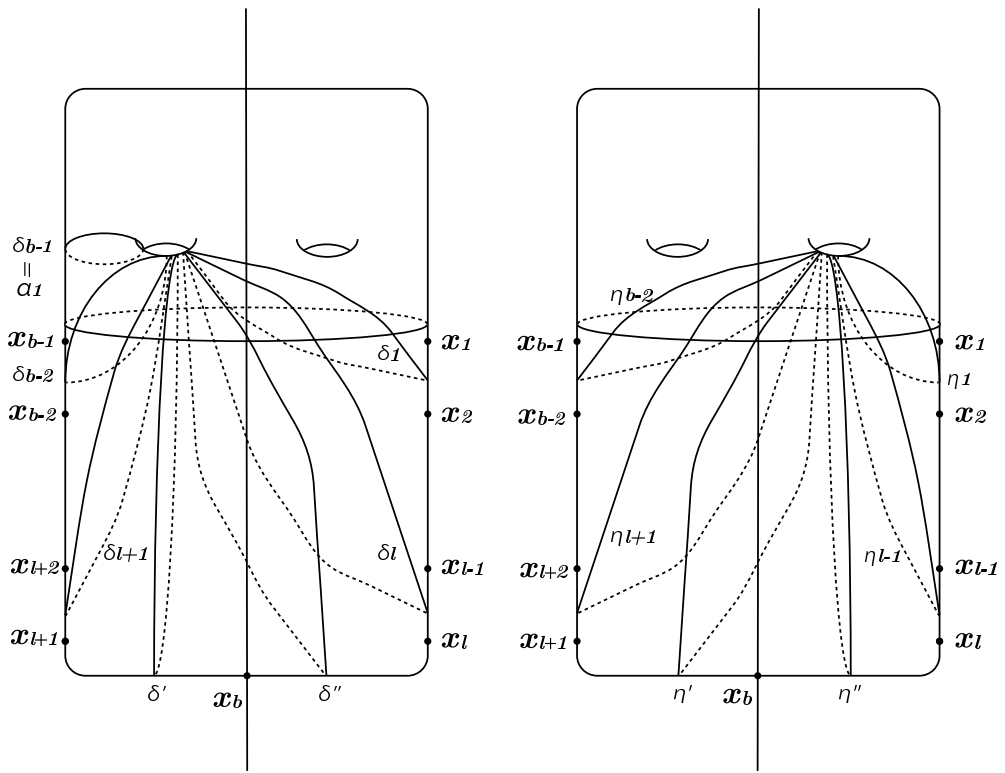}
 \end{center}
 \caption{The curves $\delta_{i}$-es, $\eta_{i}$-es.}
 \label{fig:one}
\end{figure}
\begin{lemma}
$R^{-1}(\delta_{\textit{j}})=\eta_{\textit{j}-1} \ (\textit{l}+2\leq \textit{j}\leq \textit{b}-1)$, \ $R^{-1}(\delta_{\textit{l}+1})=\eta^{\prime}$.
\end{lemma}
\begin{proof}
Figure 9 and Figure 10 shows the action of $\rho_{1}$ and $\rho_{2}$ on the curve $\delta^{\prime}$ and $\delta_{j}$ $(j=l-1,\ldots ,b-1)$. It is clear from the picture that $\eta_{j-1}=\rho_{1}\rho_{2}(\delta_{j})=R^{-1}(\delta_{j})$. It is also showed that $R^{-1}(\delta_{l+1})=\eta^{\prime}$.
\end{proof}
\begin{figure}[htbp]
 \begin{center}
  \includegraphics*[width=13cm]{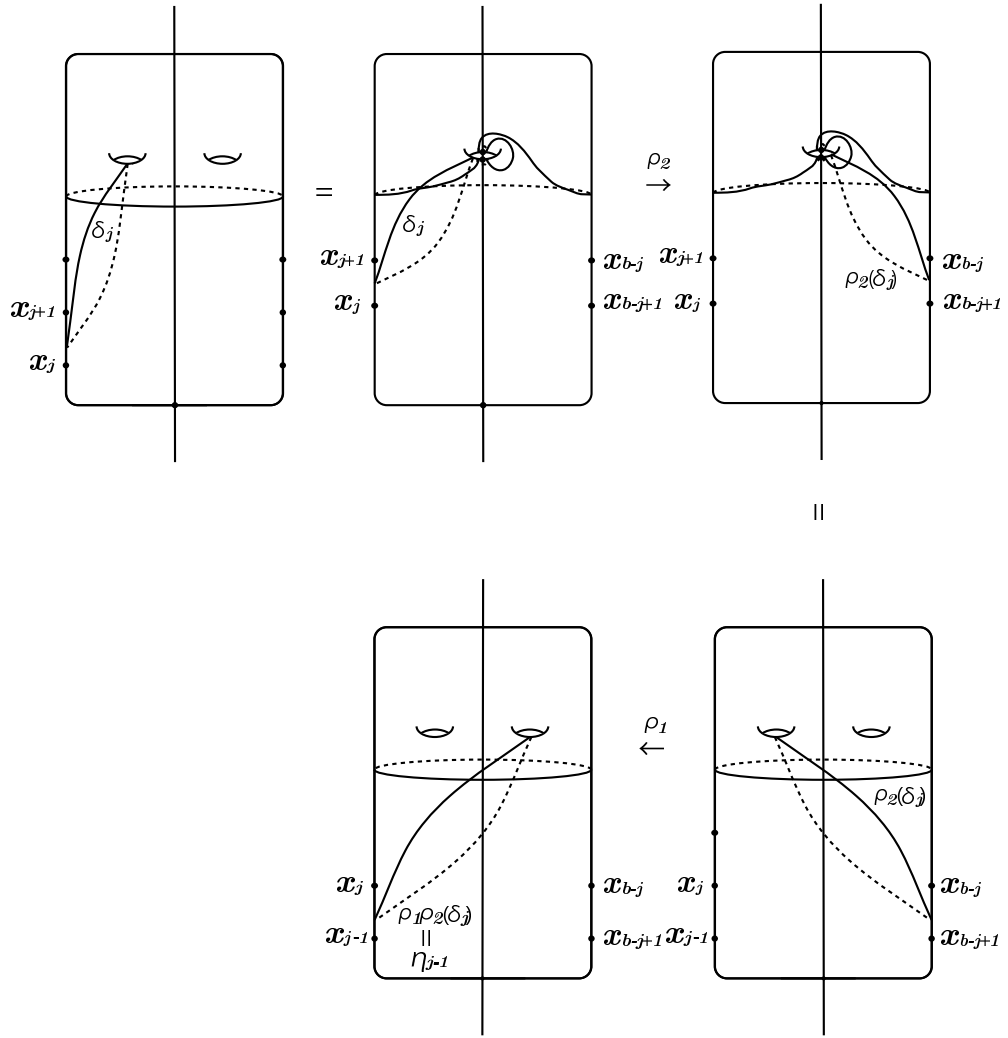}
 \end{center}
 \caption{The action of $R$ on the curve $\delta_{j}$ \ $(l+2\leq j\leq b-1)$}
 \label{fig:one}
\end{figure}
\begin{figure}[htbp]
 \begin{center}
  \includegraphics*[width=13cm]{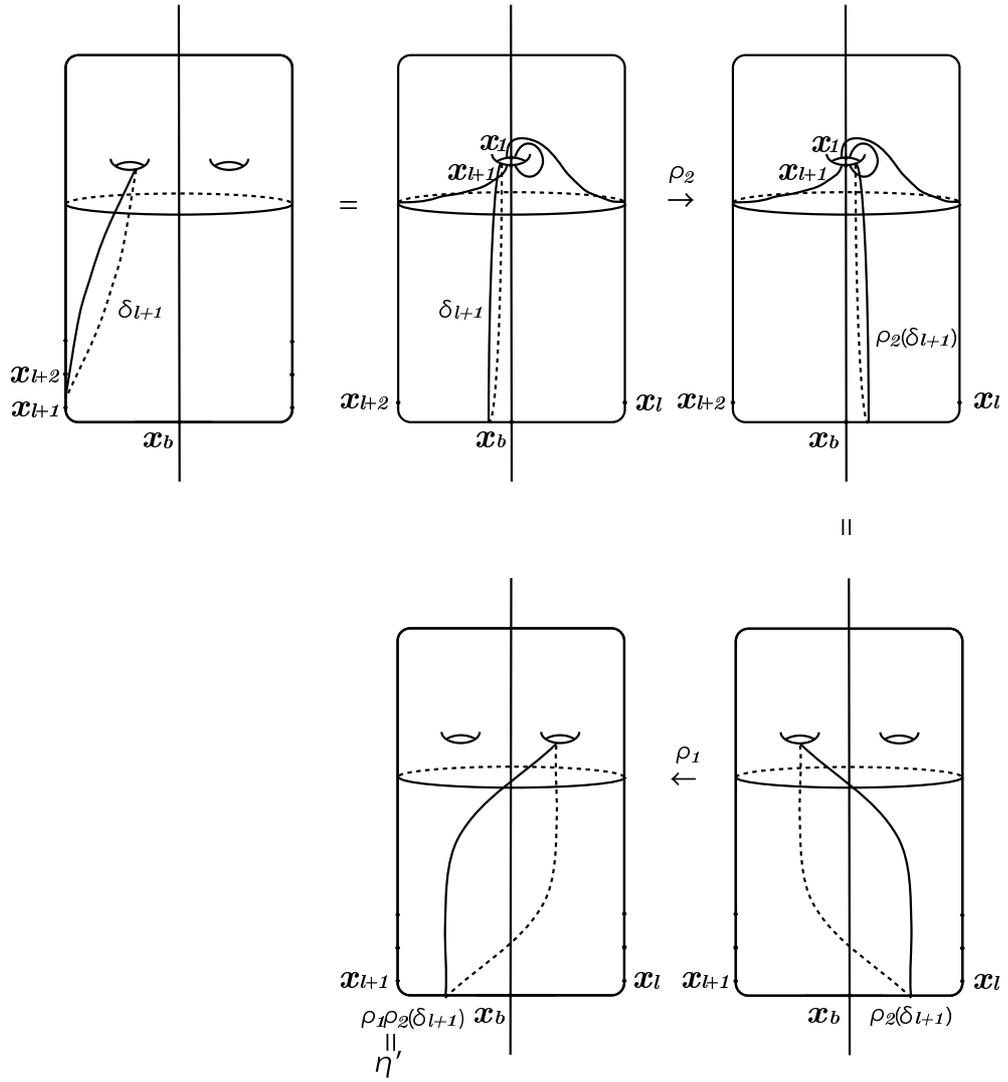}
 \end{center}
 \caption{The action of $R$ on the curve $\delta_{l+1}$} 
 \label{fig:one}
\end{figure}
\begin{figure}[htbp]
 \begin{center}
  \includegraphics*[width=13cm]{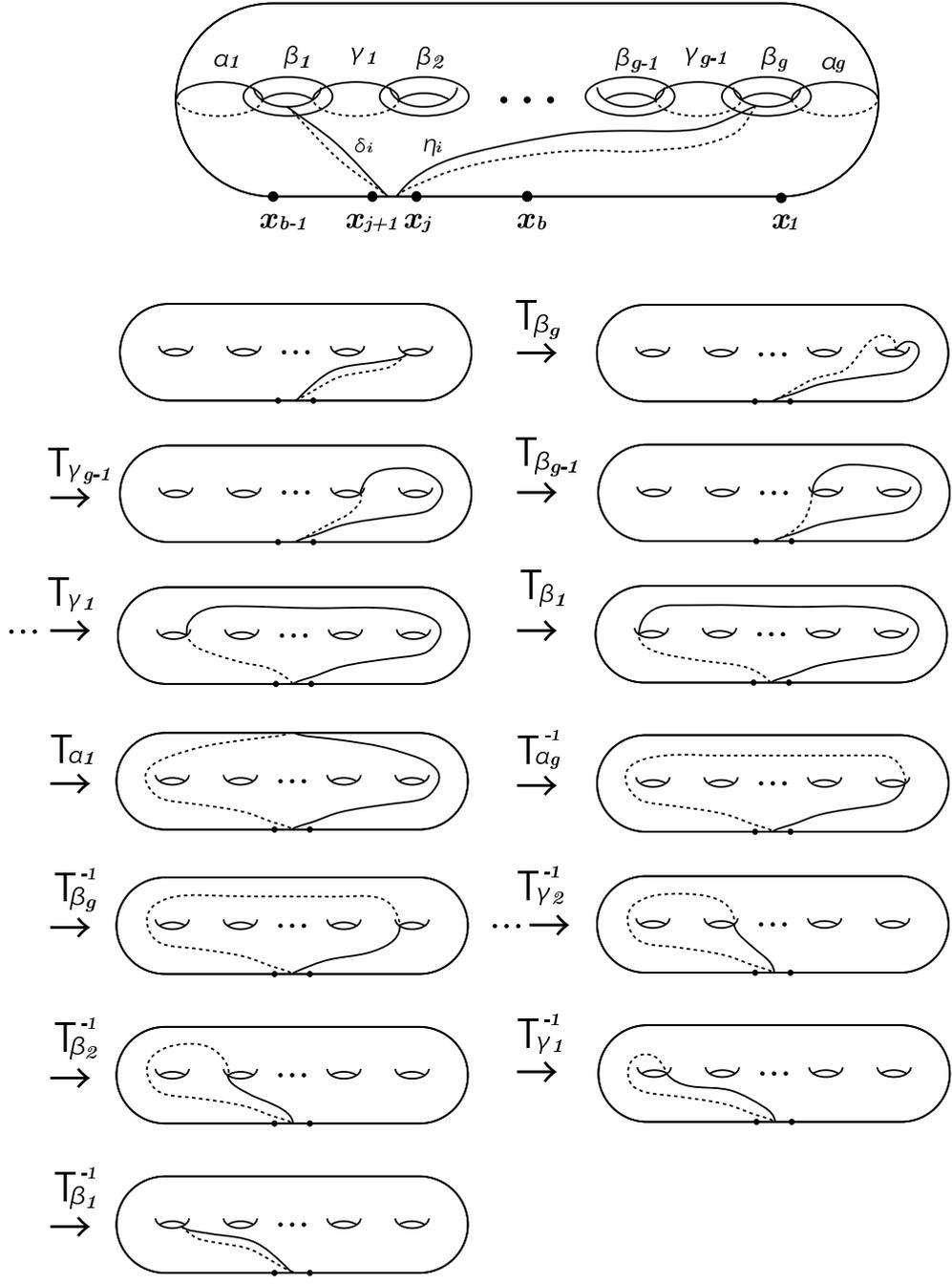}
 \end{center}
 \caption{The action of $U$ on the curve $\delta_{i}$.}
 \label{fig:one}
\end{figure}
\begin{lemma}
$T_{\delta_{\textit{j}}}, T_{\delta^{\prime}}, T_{\delta^{\prime\prime}}\in G \ \ (\textit{j}=1,\ldots ,\textit{l}-1,\textit{l}+1,\ldots,\textit{b}-2)$.
\end{lemma}
\begin{proof}
At first, we prove $T_{\delta_{j}}\in G$ using induction $j$ $(j=l+2,\ldots ,b-1)$ and $T_{\delta^{\prime}}\in G$, then we prove that $T_{j}\in G$ $(j=0,\ldots ,l-1)$ and $T_{\delta^{\prime\prime}}\in G$.\\
\ \ \ The base case, $j=b-1$, is clear because by construction $G$ contains $T_{\delta_{b-1}}=T_{\alpha_{1}}$. Suppose that $G$ contains the twist $T_{\delta_{j}}$. Using Lemma 1 and Lemma 8 we can see that the twists
$T_{\eta_{j-1}}$ and $T_{\eta^{\prime}}$ also lies in $G$, since it is conjugate to $T_{\delta_{j}}$ and $T_{\delta_{l+1}}$
\begin{eqnarray}
T_{\eta_{\textit{j}-1}}=R^{-1}T_{\delta_{\textit{j}}}R\in G.
\end{eqnarray}
Let $U$ denote the product 
\begin{eqnarray*}
U= T_{\beta_{1}}^{-1}T_{\gamma_{1}}^{-1}T_{\beta_{2}}^{-1}\cdots T_{\beta_{\textit{g}-1}}^{-1}T_{\gamma_{\textit{g}-1}}^{-1}T_{\beta_{\textit{g}}}^{-1}T_{\alpha_{\textit{g}}}^{-1}T_{\alpha_{1}}T_{\beta_{1}}T_{\gamma_{1}}T_{\beta_{2}}\cdots T_{\beta_{\textit{g}-1}}T_{\gamma_{\textit{g}-1}}T_{\beta_{\textit{g}}} \in G.
\end{eqnarray*}
The Figure 11 shows that 
\begin{eqnarray}
 U(\eta^{\prime})&=&\delta^{\prime} \nonumber \\
 U(\eta^{\prime\prime})&=&\delta^{\prime\prime}  \\
 U(\eta_{k})&=&\delta_{j} \ (j=1,\ldots,l-1,l+1,\ldots,b-2). \nonumber
\end{eqnarray}
By Lemma 1, Lemma 10 and (9) we can see that $T_{\delta_{j-1}}=UT_{\eta_{j-1}}U^{-1} \in G$. 
Thus $T_{\delta_{j}} \in G$ \ $(j=l+2,\ldots ,b-1)$. Moreover, we can see that
\begin{eqnarray*}
R^{-1}(\delta_{l+2})=\eta_{l+1},
U(\eta_{l+1})=\delta_{l+1},
R^{-1}(\delta_{l+1})=\eta^{\prime},
U(\eta^{\prime})=\delta^{\prime}.
\end{eqnarray*}
Therefore, we have that $T_{\delta_{l+1}}, T_{\delta^{\prime}}\in G$.\\
\ \ \ The next step we will prove that $T_{\delta^{\prime\prime}},T_{\delta_{\textit{j}}} \in G \ (\textit{j}=1,\ldots,\textit{l}-1)$.\\
\ \ \ By Figure 3 we can see that $\rho_{1}(\delta^{\prime\prime})=\eta^{\prime}, \rho_{1}(\delta_{\textit{j}})=\eta_{\textit{b}-1-j} \ (1 \leq \textit{j} \leq l-1)$.
By (9) we can understand that $T_{\delta_{j}}=U^{-1}T_{\eta_{b-1-j}}U, T_{\delta^{\prime\prime}}=U^{-1}T_{\eta^{\prime\prime}}U \in G$. We finished proving Lemma 10.
\end{proof}
\begin{corollary}
The group $G$ contains the subgroup $\rm Mod^{0}(\Sigma_{\textit{g,b}})$.
\end{corollary}
Therefore, we can prove Lemma 5.
\begin{remark}
The group $G$ is not the whole group $\rm Mod^{0}(\Sigma_{\textit{g,b}})$ (if $b>3$), because its image in the symmetric group $\rm Sym_{\textit{b}}$ is generated by two involutions
and therefore is a dihedral group. It is easy see that this image is the group $D_{2b}$ which is a proper subgroup of $\rm Sym_{\textit{b}}$.
\end{remark}

\section{Remark}
\ \ \ Clearly $\rm Mod(\Sigma_{\textit{g,b}})$ is never generated by two involutions, for then it would be a quotient of the infinite dihedral group, and so would be virtually abelian.
Since the current known bounds are so close to being sharp, it is natural to ask for the sharpest bounds.
\begin{problem}
For each $g\geq 3$, prove sharp bounds for the minimal number of involutions required to generate $\rm Mod(\Sigma_{\textit{g,b}})$. In particular, for $g\geq 7$ determine whether or not $\rm Mod(\Sigma_{g,b})$ is generated by 3 involutions.
\end{problem}
In order to generate the extended mapping class group $\rm Mod^{\pm}(\Sigma_{\textit{g,b}})$, it suffices to add one more generator, namely the isotopy class of any orientation-reversing diffeomorphism. Therefore, by replacing the involution $\rho_{1}$ with the reflection, we can get following result;
\begin{corollary}
For all $g\geq 3$ and $b\geq 0$, the extended mapping class group $\rm Mod^{\pm}({\Sigma_{\textit{g,b}}})$ can be generated by:\\
\ \ \ $(a)$ $4$ involutions if $g\geq 7$;\\
\ \ \ $(b)$ $5$ involutions if $g\geq 5$.
\end{corollary}
But in the case of $b=0$, Stukow showed the result that was stronger than Corollary 12. He proved that $\rm Mod^{\pm}(\Sigma_{\textit{g},0})$ is generated by three involutions. Then, we can consider following problem;
\begin{problem}
Can the extended mapping class group $\rm Mod^{\pm}(\Sigma_{\textit{g,b}})$ be generated by 3 involutions ?
\end{problem}

\section{Acknowledgement}
\ \ \ I would like to thank Professor Hisaaki Endo for careful readings and for many helpful suggestions
and comments. And I would like to thank Hitomi Fukushima, Yeonhee Jang and Kouki Masumoto for many advices.

Department of Mathematics, Graduate School of Science, Osaka University, Toyonaka, Osaka 560-0043, Japan\\
\textit{E-mail adress:} \bf{n-monden@cr.math.sci.osaka-u.ac.jp}
\end{document}